\documentclass[12pt,reqno,twoside]{amsart}
\usepackage{amsthm,amsmath,amssymb,mathrsfs,amsbsy,mathabx}
\usepackage{graphicx,epsfig,wrapfig,sidecap,subfig}
\usepackage[all]{xy}
\usepackage{proof}
\usepackage{xcolor}
\usepackage{hyperref}
\usepackage{bm}
\usepackage[capitalize]{cleveref}
\usepackage{listings}

\usepackage[a4paper, hmargin=2.4cm, vmargin={3cm, 3cm}]{geometry}

\xyoption{dvips}
\xyoption{color}

\usepackage{color}

\title[Statistics for Random Partitions and Multiscale Tilings]{Statistics and Gap Distributions in Random Kakutani Partitions and Multiscale Substitution Tilings}
\date{}
\author{Yotam Smilansky}
\address{Rutgers University, New Brunswick, NJ, USA, {\tt yotam.smilansky@rutgers.edu}}

\newcommand{\N}{{\mathbb{N}}}
\newcommand{\Z}{{\mathbb{Z}}}
\newcommand{\Q}{{\mathbb {Q}}}
\newcommand{\R}{{\mathbb{R}}}
\newcommand{\C}{{\mathbb{C}}}

\newcommand{\vol}{\mathrm{vol}}

\newcommand {\ignore}[1]  {}

\theoremstyle{plain}
\newtheorem{thm}{Theorem}[section]
\newtheorem*{thm*}{Theorem}
\newtheorem{lem}[thm]{Lemma}
\newtheorem{prop}[thm]{Proposition}
\newtheorem*{prop*}{Proposition}
\newtheorem{cor}[thm]{Corollary}

\theoremstyle{definition}

\newtheorem{example}[thm]{Example}
\newtheorem{remark}[thm]{Remark}

\newtheorem*{remark*}{Remark}

\numberwithin{equation}{section}
\swapnumbers

\newif\ifdraft\drafttrue
\usepackage{color}

\begin{document}
	
\maketitle
\begin{abstract}
		We study statistics of tiles in random incommensurable Kakutani sequences of partitions in $\mathbb{R}^d$. We provide explicit formulas that illustrate the dependence on the combinatorial structure, the volumes of the participating tiles and the entropy of the partitions in the underlying random substitution system. These improve previous results for non-random Kakutani partitions and multiscale substitution tilings, and imply a gap distribution formula for Delone sets associated with multiscale substitution tilings of the real line.  
\end{abstract}

\noindent {\bf Keywords:} tile frequencies, gap distributions, incommensurable  substitution systems, Kakutani splitting procedure,  multiscale substitution tilings, random partitions.

\section{Introduction}\label{sec: intro}

The study of the distribution of gaps in sequences on the real line is well-known to shed light on their structural properties and the extent of ``order'' they demonstrate. The same is true for the distribution of types and scales of tiles in sequences of partitions and tilings of $\R^d$, which in the latter case is related to dynamical properties such as the unique ergodicity of the associated dynamical system and spectral properties of the associated mathematical diffraction measure. In this paper we consider these statistics in sequences and tilings generated by random and non-random multiscale substitution systems, defined according to Kakutani's splitting procedure.

Before we discuss the general construction and results, let us begin with a simple $1$-dimensional motivating example. Fix $\alpha\in(0,1)$, and consider the following sequence of partitions of the unit interval $I$. First, define $\pi_1$ to be the partition of $I$ into two intervals, a red interval of length $\alpha$ and a blue one of length $1-\alpha$. Next, for any $m\in\N$, define $\pi_{m+1}$ to be the partition derived from $\pi_m$ by splitting all intervals of maximal length in $\pi_m$, each into two intervals, a red one of ratio $\alpha$ and a blue one of ratio $1-\alpha$, just as in the first step.  The splitting of intervals according to their length is sometimes known as {\it Kakutani's splitting procedure}. Indeed, choosing the left end-points of the intervals of each partition defines the classic $\alpha$-Kakutani sequence, which was shown by Kakutani \cite{Kakutani} to uniformly distribute in $I$ for any $\alpha\in(0,1)$. More generally, one can split red intervals according to  $\alpha$ as before, and blue intervals according to some, perhaps different, fixed constant $\beta\in(0,1)$. Interpreting the splitting of intervals as a substitution according to a fixed substitution rule on red and blue intervals, the above constructions can be viewed as simple examples of Kakutani sequences of partitions generated by multiscale substitution schemes as introduced in \cite{Yotam Kakutani}, to be defined shortly.

Adding a random ingredient to this process brings forward a further generalization of this construction. Consider red and blue intervals as above, and fix a red coin and constants $\alpha_1,\alpha_2\in(0,1)$. Tossing the red coin before every splitting of a red interval to determine if the partition is done according to $\alpha_1$ or $\alpha_2$ defines a random Kakutani sequence of partitions. Similarly, one can also fix a blue coin, or even a blue dice, add more intervals and colors, or extend this construction to higher dimensions.

Several questions concerning the expected statistics of these colored sequences of partitions naturally arise:
\begin{itemize}
	\item What is the asymptotic ratio of red intervals within the total number of intervals?
	\item What is the asymptotic measure of the union of all red intervals?
	\item Under an appropriate renormalization, how are the lengths of the intervals distributed? Equivalently, what is the gap distribution of the associated sequences?
\end{itemize}
In the non-random case, the first two questions were considered in \cite{Yotam Kakutani}.
It was shown that in general the limits may not exist, but under a certain irrationality assumption known as incommensurability, these two limits, which are also referred to as the {\it tile frequencies}, do exist and can be computed. Additionally, steps towards the third question in the non-random case appeared in a joint work with Yaar Solomon \cite{SmiSol} that concerns incommensurable multiscale substitution tilings. These are tilings of $\R^d$ that can be defined as certain limits of rescaled incommensurable Kakutani partitions. However, in both cases the precise way in which these frequencies and distributions are related to properties of the underlying substitution scheme remained unclear. 

 In the current paper we approach these questions in the case of random partitions in $\R^d$, and establish new formulas for the expected values of the frequencies and distributions described above.
 The results demonstrate the explicit dependence of the tile-statistics on the combinatorial, tile-volume and entropy information carried by the generating random substitution system. 
 Our results hold also for the non-random case, extending and improving the results from \cite{Yotam Kakutani} and \cite{SmiSol} mentioned above. We follow a similar approach here, taking advantage of the relations between tile frequencies and path counting results from \cite{Graphs}. In this respect, \cite{Graphs,Yotam Kakutani,SmiSol} contain a large part of the heavy lifting of our analysis.
 
We note that various non-random generalizations of $\alpha$-Kakutani sequences that can also be constructed via substitution rules have been previously studied with focus on equidistribution and discrepancy estimates, see \cite{Aistleitner Hofer,Drmota Infusino,Volcic} and the recent infinite extension introduced in \cite{Pollicott-Sewell}. Certain random generalizations have also been considered before, see in particular \cite{Pyke and van Zwet} and references within as well as \cite{BD1,BD2}, in which both the substitution rule and the splitting process itself are randomized. In the context of tilings of Euclidean space, earlier constructions include the generalized pinwheel tilings in \cite{Sadun - generalized Pinwhell} and the $1$-dimensional fusion tiling of \cite[\S A.5]{Frank-Sadun (ILC fusion)}. Unlike the partitions and tilings considered in this paper, all of these constructions consist of tiles of a single label, or type. 

Some of the tools we develop in this work, in particular those that concern the graph associated with a random substitution system, are useful in the study of random multiscale substitution tilings. These are currently studied jointly with Yaar Solomon, as part of the recent surge of interest in tilings that are defined according to some form of mixed or randomly chosen substitution rules, with examples including the various constructions considered in \cite{FSS,GM,RustSpindeler} and \cite{SchmiedingTrevino}.

\section{Preliminaries and main results}\label{sec: main results}

\subsection{Substitution systems, Kakutani partitions and the substitution semi-flow}

A {\it tile} is a labeled measurable set in $\R^d$ with positive measure denoted by $\vol$ and with boundary of measure zero. A {\it multiscale substitution scheme} $\sigma$ in $\R^d$ consists of a set of {\it prototiles} $\tau_\sigma=(T_1,\ldots,T_n)$, and for each $T_i\in\tau_\sigma$ a list of {\it substitution tiles}
\begin{equation}\label{eq: substitution tiles}
\omega(T_i)=\left(\alpha_{ij}^{(k)}T_j:\,j=1,\ldots,n,\,k=1,\ldots,k_{ij}\right)
\end{equation}
with {\it scales} $\alpha_{ij}^{(k)}\in(0,1)$ and a {\it substitution rule} $\varrho(T_i)$, which is a partition of $T_i\in\tau_\sigma$ into the elements of $\omega(T_i)$, up to isometries. One can think of $\varrho(T_i)$ as a solution to a jigsaw puzzle with pieces in $\omega(T_i)$.
A tile that arises as a copy of the prototile $T_j$, for example via a substitution, is labeled $j$ and is said to be of {\it type} $j$.  We denote by $\omega_\sigma$ and $\varrho_\sigma$ the entire collections of substitution tiles and rules in $\sigma$, respectively.

More generally, a  {\it random substitution system}  ${\bf \Sigma}$ in $\R^d$ consists of a fixed set of prototiles $\tau_{\bf \Sigma}=(T_1,\ldots,T_n)$, but now for each prototile $T_i\in\tau_{\bf \Sigma}$ there are $\ell_i\in\N$ substitution rules $\varrho_1(T_i),\ldots,\varrho_{\ell_i}(T_i)$ with substitution tiles $\omega_1(T_i),\ldots,\omega_{\ell_i}(T_i)$, and a probability vector ${\bf p}_i=(p_{i,1},\ldots,p_{i,\ell_i})$,  with $p_{i,k}>0$ and $\sum_{k=1}^{\ell_i} p_{i,k}=1$. We denote by 
$\omega_{\bf \Sigma}$ the entire collection of substitution tiles that appear in the various $\omega_k (T_i)$'s. The case in which $\ell_i=1$ for every prototile $T_i\in\tau_{\bf \Sigma}$ will be referred to as the {\it non-random case}, and then we revert to the notation for multiscale substitution schemes given above and use $\sigma$ instead of ${\bf \Sigma}$, as in \cite{Yotam Kakutani} and \cite{SmiSol}. The system ${\bf \Sigma}$ is {\it normalized} if all prototiles are of unit volume, and unless otherwise stated we will assume this throughout.

\begin{example}\label{ex: square scheme}
	A random substitution system in $\R^2$ is illustrated in Figure \ref{fig: random square scheme}. There are $n=2$ prototiles:  a white and a blue unit squares, where the white square $T_1$ has $\ell_1=2$ distinct substitution rules, and the blue square $T_2$ has $\ell_2=1$ substitution rule. 

\begin{figure}[ht!]
	\includegraphics[scale=0.9]{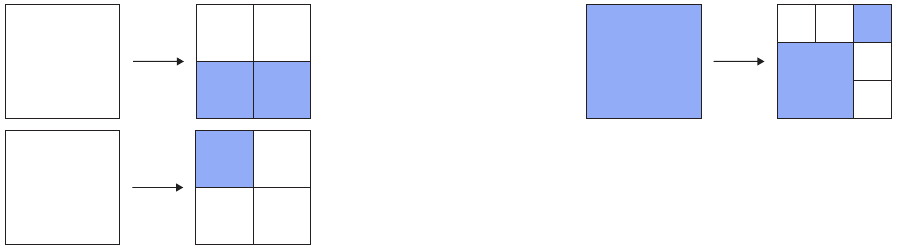}\caption{A random substitution system on two prototiles.}	\label{fig: random square scheme}
\end{figure}
\end{example}

\begin{figure}[ht!]
	\includegraphics[scale=1.04]{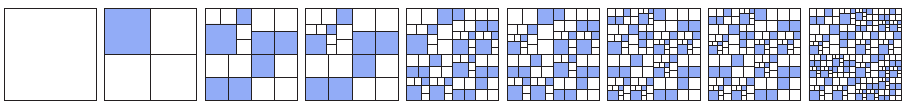}\caption{A realization of a random Kakutani sequence of partitions generated by the substitution system in Example \ref{ex: square scheme}, with $p_{1,1}=p_{1,2}=1/2$.}	\label{fig: random square Kakutani}
\end{figure}

Given a random substitution system ${\bf \Sigma}$, we define a {\it random Kakutani sequence of partitions} $\left(\pi^{{\bf \Sigma}}_m\right)_{m\ge0}$ of $T_i\in\tau_{\bf \Sigma}$. Setting $\pi^{\bf \Sigma}_0:=T_i$, for every $m\in\N$ we define $\pi^{\bf \Sigma}_{m+1}$ according to Kakutani's splitting procedure and the system ${\bf \Sigma}$, namely by substituting {\bf tiles of maximal volume} in $\pi_m$, where for each such individual tile, if it is of type $i$ a substitution rule $\varrho_k(T_i)$ is randomly chosen according to ${\bf p}_i$, independently of other tile substitutions. The first few elements in a Kakutani sequence of partitions generated by the substitution system described in Example \ref{ex: square scheme} are illustrated in Figure \ref{fig: random square Kakutani}.

It is convenient to embed the discrete procedure that defines Kakutani sequences of partitions within a continuous process that depends on a continuous parameter. Given a random substitution system ${\bf \Sigma}$, we define a {\it random substitution semi-flow} $F^{{\bf \Sigma}}_t$ and consider realizations of patches of the form $F^{{\bf \Sigma}}_t(S)$. Here $t\ge0$ is {\it time}, $S$ is a {\it tile of legal type and scale}, that is, a labeled tile  of volume not greater than the unit volume, and every realization $F^{{\bf \Sigma}}_t(S)$ is the result of the following random procedure. The tile $S$ is inflated exponentially in time, and once its volume is greater than the unit volume (which may happen at $t=0$ if $S$ is a prototile in $\tau_{\bf \Sigma}$), it is substituted according to ${\bf \Sigma}$, exactly as in the Kakutani sequence described above. As $t$ grows further the resulting patch is inflated by $e^t$, and whenever any tile surpasses the unit volume it is similarly substituted according to ${\bf \Sigma}$, independently of all other tiles. Equivalently, $F^{{\bf \Sigma}}_t(S)$ is a patch defined by repeatedly partitioning the inflated tile $e^tS$ and all subsequent tiles according to ${\bf \Sigma}$, independently, until all participating tiles are of unit volume or less.

\begin{figure}[ht!]
	\includegraphics[scale=1.04]{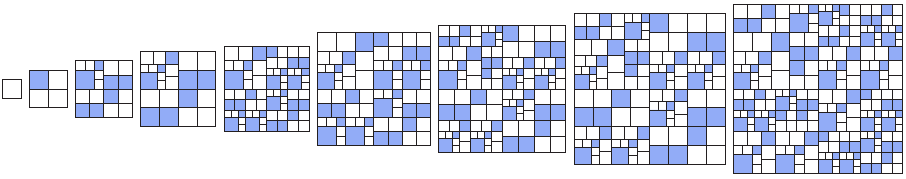}\caption{Patches in a realization of a random substitution semi-flow generated by the substitution system in Example \ref{ex: square scheme}, compare with Figure \ref{fig: random square Kakutani}.}	\label{fig: random square patches}
\end{figure}

Observe that the elements of Kakutani sequences of partitions are a countable subset of rescaled patches from a continuous collection of patches of the form $\{F^{{\bf \Sigma}}_t(T_i):\,t\ge0\}$. In fact, we have
\begin{equation}\label{eq: Kakutani = patches}
\pi^{{\bf \Sigma}}_m=\frac{1}{\vol(e^{t_m}T_i)}F_{t_m}^{{\bf \Sigma}}(T_i)
\end{equation}  
where $(t_m)_{m\ge 0}$ is the strictly increasing sequence of times for which the corresponding realizations $F_{t_m}^{{\bf \Sigma}}(T_i)$ contain a tile of unit volume. For example, patches generated by the substitution system described in Example \ref{ex: square scheme} that contain a tile of unit volume are illustrated in Figure \ref{fig: random square patches}, and are simply the partitions in Figure \ref{fig: random square Kakutani} rescaled.

The volumes of tiles in patches of the form $F^{{\bf \Sigma}}_t(S)$ are bounded from below by a positive constant depending on ${\bf \Sigma}$ and from above by $1$. It follows that the substitution semi-flow provides us with an appropriate normalization for the study of the distributions of tile volumes in Kakutani sequences of partitions. 

Our main results concern tile statistics of the above families of random patches, and applications to random and non-random Kakutani sequences of partitions, multiscale substitution tilings and their associated point processes. In order to state the main results we must first define the graphs and matrices that carry the relevant information of the underlying substitution system.

\subsection{Graphs and matrices associated with random substitution systems}

The directed graph $G_\sigma$ associated with a multiscale substitution scheme  $\sigma$ has vertices $\{1,\ldots,n\}$ corresponding to prototiles $\tau_\sigma=(T_1,\ldots,T_n)$, and edges defined by the substitution tiles $\omega_{\sigma}$, where every substitution tile $T=\alpha T_j\in \omega(T_i)$ defines a directed edge $\varepsilon_T$ with initial vertex $i$, terminal vertex $j$ and length $l(\varepsilon_T)=\log \frac{1}{\vol T}$. Similarly, the graph $G_{\bf \Sigma}$ associated with a random substitution system  ${\bf \Sigma}$ has vertices corresponding to the prototiles $\tau$, and edges defined by all substitution tiles $\omega_{\bf \Sigma}$ combined.

If the graph $G_\sigma$ (respectively, $G_{\bf \Sigma}$) is strongly connected, then the substitution scheme $\sigma$ (system ${\bf \Sigma}$) is said to be {\it irreducible}. If it contains two orbits of incommensurable lengths $a,b$ with $a\notin b\mathbb{Q}$, then the graph and the scheme (system) are both said to be {\it incommensurable}. Unless otherwise stated we assume from here on that ${\bf \Sigma}$ is {\bf normalized, irreducible and incommensurable}.  

\begin{remark}\label{rem: different definition of G}
	For a detailed discussion and many illustrated examples concerning multiscale substitution schemes, their basic properties and their associated graphs we refer the reader to \cite{Yotam Kakutani,SmiSol}. Note that our definition of $G_\sigma$ in the current paper is slightly different from the definition used previously, where the length of $\varepsilon_T$ was defined to be $\log\frac{1}{\alpha}=\frac{1}{d}\log\frac{1}{\vol T}$. The two definitions give rise to graphs with identical combinatorial structure, and with the edges of one simply rescaled copies of the edges of the other by a factor of the dimension $d$. We work with one definition and not the other only to simplify our presentation. 
\end{remark}

Given a substitution scheme $\sigma$, the {\it substitution matrix} $S_\sigma\in M_n(\Z)$ is defined by 
		\begin{equation*}\label{eq: substitution matrix}
		\left(S_\sigma\right)_{ij}=\#\left\{\text{Rescaled copies of }T_j \text{ in } \omega(T_i)\right\}=\sum\limits_{\substack{T\in\omega(T_i) \\ T\text { of type }j}}1,
		\end{equation*} 
		 and carries the combinatorial information of the substitution scheme. Note that $S_\sigma$ is the adjacency matrix of $G_\sigma$ when considered as a combinatorial object. The {\it volume matrix} $V_\sigma\in M_n(\R)$ is defined by
		\begin{equation*}\label{eq: volume matrix}
		\left(V_\sigma\right)_{ij}=\sum\limits_{\substack{T\in\omega(T_i) \\ T\text { of type }j}}\vol T.
		\end{equation*} 
		Since $\sigma$ is normalized ${\bf u}_\vol={\bf 1}:=(1,\ldots,1)^T\in\R^n$ is a right Perron--Frobenius eigenvector of $V_\sigma$, with eigenvalue $\mu=1$. The {\it entropy matrix} $H_\sigma\in M_n(\R)$ is defined by
		\begin{equation*}\label{eq: entroy matrix}
		\left(H_\sigma\right)_{ij}=\sum\limits_{\substack{T\in\omega(T_i) \\ T\text { of type }j}}-\vol T\cdot\log\vol T. 
		\end{equation*} 
		Note that the entries of $H_\sigma{\bf 1}$, the vector of sums of rows of $H_\sigma$, are exactly the partition entropies of the partitions $\varrho(T_i)$ of $T_i\in\tau_\sigma$ for $1\le i\le n$. Therefore $\left(H_\sigma\right)_{ij}$ can be thought of as the contribution of tiles of type $j$ to the entropy of the partition $\varrho(T_i)$.
		
		The substitution, volume and entropy matrices associated with a random substitution system ${\bf \Sigma}$ are similarly defined, weighing the contributions of each of the substitution rules on each prototile according to the corresponding probability vectors: 
		\begin{align}\label{eq: S,V,H random}
		\left(S_{{\bf \Sigma}}\right)_{ij}&=\sum_{k=1}^{\ell_i} p_{i,k}\sum\limits_{\substack{T\in\omega_k(T_i) \\ T\text { of type }j}}1\qquad
		\left(V_{{\bf \Sigma}}\right)_{ij}=\sum_{k=1}^{\ell_i} p_{i,k}\sum\limits_{\substack{T\in\omega_k(T_i) \\ T\text { of type }j}}\vol T\\
		\left(H_{{\bf \Sigma}}\right)_{ij}&=\sum_{k=1}^{\ell_i} p_{i,k}\sum\limits_{\substack{T\in\omega_k(T_i) \\ T\text { of type }j}}-\vol T\cdot \log \vol T.\notag
		\end{align}  		
		
In addition, we define the matrix valued functions $C_{{\bf \Sigma}},D_{{\bf \Sigma}}:\R\rightarrow M_n(\R)$ by
		\begin{equation}\label{eq: density matrices}
		\left(C_{{\bf \Sigma}}(x)\right)_{ij}=\sum_{k=1}^{\ell_i} p_{i,k}\sum\limits_{\substack{T\in\omega_k(T_i) \\ T\text { of type }j}}c_{\vol T}(x) \qquad \left(D_{{\bf \Sigma}}(x)\right)_{ij}=\sum_{k=1}^{\ell_i} p_{i,k}\sum\limits_{\substack{T\in\omega_k(T_i) \\ T\text { of type }j}}d_{\vol T}(x),
		\end{equation}
		where
		\begin{equation}\label{eq: density functions}
		c_{\vol T}(x)=\begin{cases}
		\frac{\vol T}{x^2},&\vol T<x\le 1  \\
		0,&\rm{otherwise} 
		\end{cases} \qquad d_{\vol T}(x)=\begin{cases}
		\frac{\vol T}{x},&\vol T<x\le 1  \\
		0,&\rm{otherwise} 
		\end{cases}.
		\end{equation}
		Note that the entries of $C_{{\bf \Sigma}}(x)$ and $D_{{\bf \Sigma}}(x)$ are piecewise smooth functions with support bounded away from zero. 

\subsection{Main results}
The following formulas constitute our main result.
\begin{thm}\label{thm: Random formulas}
	Let ${\bf \Sigma}$ be a normalized irreducible incommensurable random substitution system in $\R^d$. Let  $S$ be a tile of legal type and scale, let $1\le r\le n$ and let $0\le a\le b\le 1$. Then
	\begin{align}
	\frac{\mathbb{E}\left[\#\{\emph{Tiles of type } r \emph{ and }\vol\in[a,b]\emph{ in } F^{{\bf \Sigma}}_t(S)\}\right]}{\vol(e^tS)}&=\int_a^b\frac{\left[{\bf v}^TC_{\bf \Sigma}(x)\right]_r}{{\bf v}^TH_{\bf \Sigma} {\bf 1}}dx+o(1),\,t\rightarrow\infty,\notag\\
	\frac{\mathbb{E}\left[\vol\left(\bigcup\emph{Tiles of type } r \emph{ and }\vol\in[a,b]\emph{ in } F^{{\bf \Sigma}}_t(S)\right)\right]}{\vol(e^tS)}&=\int_a^b\frac{\left[{\bf v}^TD_{\bf \Sigma}(x)\right]_r}{{\bf v}^TH_{\bf \Sigma} {\bf 1}}dx+o(1),\,t\rightarrow\infty,\notag
	\end{align}
	independently of the type of $S$, where ${\bf v}\in \R^n$ is a left Perron--Frobenius eigenvector of $V_{{\bf \Sigma}}$ and ${\bf 1}=(1,\ldots,1)^T\in\R^n$. In particular
	\begin{align*}
	\frac{\mathbb{E}\left[\#\{\emph{Tiles of type } r \emph{ in } F^{{\bf \Sigma}}_t(S)\}\right]}{\vol(e^tS)}&=\frac{\left[{\bf v}^T\left(S_{\bf \Sigma}-V_{\bf \Sigma}\right)\right]_r}{{\bf v}^TH_{\bf \Sigma} {\bf 1}}+o(1),\,t\rightarrow\infty,\\
	\frac{\mathbb{E}\left[\vol\left(\bigcup\emph{Type } r \emph{ in } F^{{\bf \Sigma}}_t(S)\right)\right]}{\vol(e^tS)}&=\frac{\left[{\bf v}^TH_{\bf \Sigma}\right]_r}{{\bf v}^TH_{\bf \Sigma} {\bf 1}}+o(1),\,t\rightarrow\infty.
	\end{align*}
	In the non-random case the expected value symbols are redundant.
\end{thm}

Note that the formulas in Theorem \ref{thm: Random formulas} do not depend on the geometric nature of the tiles or on the configuration of tiles in the substitution rules $\varrho_k(T_i)$, only on the information recorded by the matrices associated with ${\bf \Sigma}$, namely the number of tiles of each type in each $\omega_k(T_i)$, on their volumes and on the partition entropy, and therefore hold in a much more general setting. 

 
 The following Corollary \ref{cor: random kakutani frequencies} is immediate in view of \eqref{eq: Kakutani = patches}.

\begin{cor}\label{cor: random kakutani frequencies}
	Let ${\bf \Sigma}$ be a normalized irreducible incommensurable random substitution system in $\R^d$. Let  $\left(\pi^{\bf \Sigma}_m\right)_{m\ge0}$ be a random Kakutani sequence of partitions of $T_i\in\tau_\sigma$, and let $1\le r\le n$. Then 
	\begin{align*}
	\frac{\mathbb{E}\left[\#\{\emph{Tiles of type } r \emph{ in } \pi^{\bf \Sigma}_m\}\right]}{\mathbb{E}\left[\#\{\emph{Tiles in } \pi^{\bf \Sigma}_m\}\right]}&=\frac{\left[{\bf v}^T\left(S_{\bf \Sigma}-V_{\bf \Sigma}\right)\right]_r}{{\bf v}^T\left(S_{\bf \Sigma}-V_{\bf \Sigma}\right) {\bf 1}}+o(1),\\
	\mathbb{E}\left[\vol\left(\bigcup\emph{Type } r \emph{ in } \pi^{\bf \Sigma}_m\right)\right]&=\frac{\left[{\bf v}^TH_{\bf \Sigma}\right]_r}{{\bf v}^TH_{\bf \Sigma} {\bf 1}}+o(1).
	\end{align*}
	In the non-random case the expected value symbols are redundant. 
\end{cor}

We note that the existing tile frequency formulas that appear \cite[\S7]{SmiSol} and \cite[Theorem 1.13]{Yotam Kakutani} are challenging to decipher and give little additional insight to the precise way the frequencies depend on the generating substitution scheme and its basic properties. While allowing for computations in specific cases and examples, the formulas remained somewhat opaque, and were not given explicitly in terms of the associated matrices $S_\sigma, V_\sigma$ and $H_\sigma$ or in terms of the density functions $C_\sigma(x)$ and $D_\sigma(x)$, introduced above. Therefore, even in the non-random case our current contribution improves and illuminates previous results.

\subsubsection{Multiscale substitution tilings}

Tilings of $\R^d$ are {\it generated} by a substitution scheme $\sigma$ as limits of patches of the form $F^{\sigma}_t(S)$, taken with respect to a suitable topology. If we assume that the tiles in the substitution rules $\varrho_\sigma$ are all translated copies of the substitution tiles $\omega_\sigma$, then these tilings consist of tiles that are rescaled translated copies of the prototiles $\tau_\sigma$ and appear in a bounded set of scales. Moreover, tiles appear in uniform frequencies that are independent of the tiling itself and of the {\it van Hove sequence} of sets in $\R^d$ used to approach the limit, where a sequence of sets  $(A_q)_{q\ge 1}$ is van Hove if the volumes of the $A_q$'s grow faster than the volumes of the $r$-neighborhoods of their boundaries, for any $r>0$. A thorough introduction to this construction and proofs of these results appear in \cite{SmiSol}. In view of our main result, Corollary \ref{cor: Tilings formulas}  provides explicit formulas for these frequencies.

\begin{cor}\label{cor: Tilings formulas}
	Let ${\sigma}$ be a normalized irreducible incommensurable substitution scheme in $\R^d$, so that all tiles in  $\varrho_\sigma$ are translated copies of the substitution tiles $\omega_\sigma$.  Let  $\mathcal{T}$ be a multiscale substitution tiling generated by $\sigma$, let $S$ be a tile of legal type and scale, let $1\le r\le n$ and let $0\le a<b\le 1$. Then for every van Hove sequence $(A_q)_{q\ge 1}$ in $\R^d$ 
	\begin{align}
		\frac{\#\{\emph{Tiles of } \mathcal{T} \emph{ of type } r \emph{ and }\vol\in[a,b]\emph{ contained in } A_q+h\}}{\vol(A_q)}&=\int_a^b\frac{\left[{\bf v}^TC_\sigma(x)\right]_r}{{\bf v}^TH_\sigma {\bf 1}}dx+o(1),\notag\\
	\frac{\vol\left(\bigcup\emph{Tiles of } \mathcal{T} \emph{ of type } r \emph{ and }\vol\in[a,b]\emph{ contained in } A_q+h\right)}{\vol(A_q)}&=\int_a^b\frac{\left[{\bf v}^TD_\sigma(x)\right]_r}{{\bf v}^TH_\sigma {\bf 1}}dx+o(1),\notag
		\end{align}
	uniformly in $h\in\R^d$, independently of the type of $S$. 
\end{cor}

The formulas in \cite[Theorem 8 parts (4),(5)]{Sadun - generalized Pinwhell} and \cite[Theorem A.3,  for $n=1$]{Frank-Sadun (ILC fusion)} can be recovered using the above results. Note that since the aforementioned examples consist of a single prototile, all the matrices involved are of order $1$.

\begin{remark}
	The error terms in the formulas of Theorem \ref{thm: Random formulas} are large. In fact, they are at least of order $1/\log^k\vol(e^tS)$ for some fixed $k\ge1$, which depends only on the underlying system. See \cite[\S8]{SmiSol} for more details, and for a discussion about how these ``bad'' error terms are related to questions on bounded displacement and bilipschitz equivalence classes of the associated Delone sets. 
\end{remark}

\subsubsection{Gap distributions} 

In the study of ordered sequences of real numbers, two types of distributions of pairs of points play an important role. The {\it gap distribution}, which is the distribution of the distance between nearest neighbors, and the {\it pair correlation}, the distribution of distances between pairs of points in general. These distributions give information on how ``random'' or ``ordered'' a sequence is, see \cite{Marklof} for an introductory discussion with many examples. In particular they are of special interest in the study of aperiodic order, in which an important tool in measuring the extent to which a system is ``ordered'' is given by the closely related {\it mathematical diffraction}, see \cite[\S9]{Baake-Grimm} for a detailed introduction. While the gap distribution and pair correlation are of course related, it is often impossible to deduce one from the other, and in different contexts one of the two is often more accessible.

Consider an incommensurable substitution scheme $\sigma$ in $\R$ in which all participating tiles are intervals, perhaps of $n$ distinct types. Then all tiles in multiscale substitution tilings generated by $\sigma$ are intervals, and the boundary points of these intervals define Delone sets for which the following explicit gap distribution formula holds. 

\begin{cor}\label{cor: gap distributions}
	Let $\mathcal{P}$ be a Delone set associated with a normalized irreducible incommensurable substitution scheme $\sigma$ in $\R$. Then
	\begin{align*}
	\frac{\#\{\emph{Neighbors in }\mathcal{P}\cap[-N,N]\emph{ of distance in }[a,b]\}}{\#\{\mathcal{P}\cap[-N,N]\}}\rightarrow\int_a^b\frac{{\bf v}^TC_\sigma(x)\bf{1}}{{\bf v}^TH_\sigma\bf{1}}dx+o(1),\,N\rightarrow\infty.
	\end{align*}  
\end{cor}

For example, let $\alpha\in(0,1)$ with $\log\alpha\notin\log(1-\alpha)\Q$, and assume without loss of generality that $\alpha\le1-\alpha$. The density function for the gap distribution in the Delone set associated with the $\alpha$-Kakutani sequence described in \S\ref{sec: intro} is given by 
\begin{align*}
	f(x)=\begin{cases}
	\frac{1}{-\alpha\log\alpha-(1-\alpha)\log(1-\alpha)}\cdot \frac{\alpha}{x^2},&\alpha<x\le1-\alpha  \\
	\frac{1}{-\alpha\log\alpha-(1-\alpha)\log(1-\alpha)}\cdot \frac{1}{x^2},&1-\alpha<x\le1  \\
	0,&\rm{otherwise}. 
	\end{cases}
\end{align*}
Compare with \cite[Theorem A.3]{Frank-Sadun (ILC fusion)} for the case $\alpha=1/3$.

\begin{remark}
	Numerical examination suggests that unlike the case of the gap distribution, the pair correlation of such Delone sets is reminiscent of what one would expect of a random Poisson process. We hope to provide further insight on this behavior and its implications to mathematical diffraction  in future work. 
\end{remark}

\section{Counting walks in $G_{\bf \Sigma}$}\label{sec: graphs}

 Recall that substitution systems are associated with graphs. Our proof of Theorem \ref{thm: Random formulas} relies on asymptotic path counting results for incommensurable graphs, equipped with appropriate weights, established jointly with Avner Kiro and Uzy Smilansky in \cite{Graphs}.

\subsection{Counting walks in general incommensurable graphs}
Consider a graph $G$ with vertices $\{1,\ldots,n\}$ and finitely many directed edges. Denote by $w(\varepsilon)$ the weight assigned to the edge $\varepsilon$, and by $l(\varepsilon)$ its length. A {\it path} in $G$ is a directed walk that initiates and terminates at vertices of the graph. The weight of a path $\gamma$, defined to be the {\bf product of the weights} assigned to the edges that make up $\gamma$, is denoted by $w(\gamma)$, and its length, the {\bf sum of the lengths} of its edges, by $l(\gamma)$. The associated {\it graph matrix function} is the matrix valued function $M:\C\rightarrow M_n(\C)$ with entries defined by 
\begin{equation}\label{eq: graph matrix function}
\left(M(s)\right)_{ij}=\sum\limits_{\substack{\varepsilon \text{ edge}\\i\rightarrow j}}w(\varepsilon)e^{-l(\varepsilon)s}.
\end{equation}

\begin{thm}[cf. \cite{Graphs}]\label{thm: General graph counting}
	Let $G$ be a strongly connected incommensurable directed graph $G$ with weights assigned to its edges, and denote by $\Gamma_{ih}$ the set of paths  in $G$ with initial vertex $i$ and terminal vertex $h$. Let $\varepsilon$ be an edge with initial vertex $h$, and let $A$ be an interval contained in $\varepsilon$ and with boundary points at distance $\delta$ and $\eta$ of vertex $h$, for $0\le\delta\le\eta\le l(\varepsilon)$. Then the asymptotic behavior of the sum of the weights of walks of length exactly $t$ that initiate at vertex $i$ and terminate at some point on the interval $A$ is
	\begin{equation*}
\sum_{\gamma\in\Gamma_{ih}}w(\gamma)w(\varepsilon)\chi_{(l(\gamma)+\delta,l(\gamma)+\eta]}(t)= w(\varepsilon)(Q)_{ih}\int_{\delta}^{\eta}e^{-\lambda x}dx\cdot e^{\lambda t}+o\left(e^{\lambda t}\right),\,t\rightarrow\infty,
	\end{equation*}
	where  $\lambda$ is the smallest real number for which the spectral radius of $M(\lambda)$ is exactly $1$ and
	\begin{equation}\label{eq: Q adj formula}
	Q=\frac{\emph{adj}\left(I-M(\lambda)\right)}{-\emph{tr}\left(\emph{adj}\left(I-M(\lambda)\right)\cdot M'(\lambda)\right)}\in M_n(\R),
	\end{equation}
	with derivatives taken entry-wise.
\end{thm}

The proof follows from a careful study of the Laplace transform of the counting function and an application of the Wiener--Ikehara Tauberian theorem, combined with tools from the theory of Perron--Frobenius. In particular, an important step is to show that the location of the poles of the Laplace transform allows the application of the Wiener--Ikehara theorem, a property that is closely linked with the incommensurability of the graph. Full details can be found in \cite{Graphs} for the case that the interval $A$ is an entire edge, and the required background concerning the Perron--Frobenius and the Wiener--Ikehara theorems can be found in \cite[Chapter XIII]{Gantmacher} and \cite[Chapter 8.3]{Montgomery-Vaughan}, respectively. The refinement to general intervals is then simple, and can be found in the proof of \cite[Theorem 7.2]{SmiSol}. We note that though the result is stated for a half-closed half-open interval $A$, it is the same whether the boundaries of $A$ are included or not.


\subsection{Graphs associated with random substitution systems}

In the case of a graph associated with a random substitution system ${\bf \Sigma}$, the following Lemma \ref{lem: S,V,H,Q for random} allows the interpretation of the matrix $Q$ in \eqref{eq: Q adj formula} in terms of the system ${\bf \Sigma}$. 

\begin{lem}\label{lem: S,V,H,Q for random}
	Let ${\bf \Sigma}$ be a random substitution system in $\R^d$, and let $S_{\bf \Sigma}$, $V_{\bf \Sigma}$ and $H_{\bf \Sigma}$ be the associated substitution, volume and entropy matrices as in \eqref{eq: S,V,H random}. Let $G_{\bf \Sigma}$ be the associated graph, with weights assigned so that $w(\varepsilon)=p_{i,k}$ for all edges $\varepsilon$ associated with substitution tiles in $\omega_k(T_i)$. Let $M_{\bf \Sigma}(s)$ be the graph matrix function as in \eqref{eq: graph matrix function} and $Q_{\bf \Sigma}$  the corresponding matrix defined by \eqref{eq: Q adj formula}. Then  
	\begin{equation*}
	S_{\bf \Sigma}=M_{\bf \Sigma}(0),\quad V_{\bf \Sigma}=M_{\bf \Sigma}(1),\quad H_{\bf \Sigma}=-M'_{\bf \Sigma}(1), \quad\text{and}\quad Q_{\bf \Sigma}=\frac{{\bf 1} {\bf v}^T}{{\bf v}^TH_{\bf \Sigma} {\bf 1}},
	\end{equation*}
	with derivatives taken entry-wise, and where  ${\bf v}\in \R^n$ is a left Perron--Frobenius eigenvector of $V_{{\bf \Sigma}}$. In particular, the rows of $Q_{\bf \Sigma}$ are all equal. 
\end{lem}

\begin{proof}
	By the definition of $G_{\bf \Sigma}$ and the weights assigned to its edges
	\begin{equation*}
	(M_{\bf \Sigma}(s))_{ij}=\sum\limits_{\substack{\varepsilon \text{ edge}\\i\rightarrow j}}w(\varepsilon)e^{-l(\varepsilon)s}=\sum_{k=1}^{\ell_i} p_{i,k}\sum\limits_{\substack{T\in\omega_k(T_i) \\ T\text { of type }j}}(\vol T)^{s}.
	\end{equation*}
	The entries of the substitution, volume and entropy matrices can therefore be expressed in terms of the graph matrix function as
	\begin{align*}
	(S_{\bf \Sigma})_{ij}&=\sum_{k=1}^{\ell_i} p_{i,k}\sum\limits_{\substack{T\in\omega_k(T_i) \\ T\text { of type }j}}1=(M_{\bf \Sigma}(0))_{ij}\quad\qquad
	(V_{\bf \Sigma})_{ij}=\sum_{k=1}^{\ell_i} p_{i,k}\sum\limits_{\substack{T\in\omega_k(T_i) \\ T\text { of type }j}}\vol T=(M_{\bf \Sigma}(1))_{ij}\\(H_{\bf \Sigma})_{ij}&=\sum_{k=1}^{\ell_i} p_{i,k}\sum\limits_{\substack{T\in\omega_k(T_i) \\ T\text { of type }j}}-\vol T\cdot\log\vol T=- (M'_{\bf \Sigma}(1))_{ij},
	\end{align*}
	satisfying the assertions concerning $S_{\bf \Sigma}$, $V_{\bf \Sigma}$ and $H_{\bf \Sigma}$. 
	
	We now consider $Q_{\bf \Sigma}$. Since ${\bf \Sigma}$ is normalized, the volumes of all tiles in $\omega_k(T_i)$ sum up to $1$ for every prototile $T_i\in\tau_{\bf \Sigma}$ and every $1\le k\le \ell_i$. It follows that for every $1\le i\le n$
	\begin{equation*}
	\left[M_{\bf \Sigma}(1){\bf 1}\right]_i=\left[V_{\bf \Sigma}{\bf 1}\right]_i=\sum_{j=1}^n\sum_{k=1}^{\ell_i} p_{i,k}\sum\limits_{\substack{T\in\omega_k(T_i) \\ T\text { of type }j}}\vol T=\sum_{k=1}^{\ell_i} p_{i,k}\sum_{j=1}^n\sum\limits_{\substack{T\in\omega_k(T_i) \\ T\text { of type }j}}\vol T=\sum_{k=1}^{\ell_i} p_{i,k}=1.
	\end{equation*}
	That is, $M_{\bf \Sigma}(1){\bf 1}={\bf 1}$ and $M_{\bf \Sigma}(1)$ has a Perron--Frobenius eigenvalue $\mu=1$ with ${\bf u}_\vol={\bf 1}$ a right Perron--Frobenius eigenvector. In addition, clearly for any $x>1$ the sums of the rows of $M_{\bf \Sigma}(x)$ are all strictly smaller than $1$, and so the spectral radius of $M_{\bf \Sigma}(x)$ is strictly smaller than $1$, see for example \cite[p. 63]{Gantmacher}. It follows that  $\lambda=1$ in the case of a graph $G_{\bf \Sigma}$ associated with a substitution system ${\bf \Sigma}$ and equipped with weights according to the probabilities given by the ${\bf p}_i$'s.

	By the theory of Perron--Frobenius, the matrix $V_{\bf \Sigma}$ has a non-negative left Perron--Frobenius eigenvector ${\bf v}\in\R^n$ for which 
	\begin{equation}\label{eq: adjugate}
	\text{adj}\left(I-V_{\bf \Sigma}\right)={\bf 1}{\bf v}^T.
	\end{equation}
	It follows that 
	\begin{equation}\label{eq: trace adjugate}
	-\text{tr}\left(\text{adj}\left(I-V_{\bf \Sigma}\right)\cdot M_{\bf \Sigma}'(1)\right)=-\text{tr}\left(-{\bf 1}{\bf v}^TH_{\bf \Sigma}\right)=\text{tr}\left({\bf v}^TH_{\bf \Sigma}{\bf 1}\right)={\bf v}^TH_{\bf \Sigma}{\bf 1},
	\end{equation}
	where the fact that $\text{tr}(AB)=\text{tr}(BA)$ is used. Combining \eqref{eq: adjugate} and \eqref{eq: trace adjugate} we get
	\begin{equation*}\label{eq: formula for Q_sigma}
	Q_{\bf \Sigma}=\frac{\text{adj}\left(I-M_{\bf \Sigma}(1)\right)}{-\text{tr}\left(\text{adj}\left(I-M_{\bf \Sigma}(1)\right)\cdot M'_{\bf \Sigma}(1)\right)}=\frac{{\bf 1} {\bf v}^T}{{\bf v}^TH_{\bf \Sigma} {\bf 1}}.
	\end{equation*}
	 Note that this expression is independent of the choice of a left Perron--Frobenius eigenvector ${\bf v}$ of $V_{\bf \Sigma}$. Finally, since the columns of $Q_{\bf \Sigma}$ are spanned by ${\bf 1}$, its rows are equal. 
\end{proof}

\section{Proof of main result}\label{sec: proofs}

In the non-random case, if $\sigma$ is a substitution scheme in $\R^d$ then an important observation is the following correspondence between tiles in a patch of the form $F^\sigma_t(T_i)$ and walks on the associated graph: for every $t>0$, the tiles of $F^\sigma_t(T_i)$ are in one-to-one correspondence with walks of length $d\cdot t$ originating at the vertex $i$ in $G_\sigma$. Moreover, if $T$ is a tile of type $r$ and volume $\vol T$, then the corresponding walk $\gamma_T$ ends at a point on an edge that terminates at vertex $r$, and the ending point is at distance $\log\frac{1}{\vol T}$ from the vertex $r$. See \cite[\S2]{SmiSol} for additional details, examples and illustrations, and notice the slight difference in the definition of the associated graph pointed out in Remark \ref{rem: different definition of G}. 

Given a random substitution system ${\bf \Sigma}$ in $\R^d$, this correspondence extends to tiles in realizations of $F_t^{\bf \Sigma}(T_i)$ and walks of length $d\cdot t$ originating at the vertex $i$ in the associated graph $G_{\bf \Sigma}$. Every tile in a patch of the form $F_t^{\bf \Sigma}(T_i)$ corresponds to a walk of length $d\cdot t$ on $G_{\bf \Sigma}$, and every walk on $G_{\bf \Sigma}$ corresponds to a tile that may appear after applying the random substitution semi-flow on $T_i$ for time $t$, that is, a tile in a realization of $F_t^{\bf \Sigma}(T_i)$. The volume and the type of the tile corresponds to the termination point of the walk, exactly as in the non-random case. 
An appropriate choice of weights for the edges of $G_{\bf \Sigma}$ and the application of Theorem \ref{thm: General graph counting} will yield the formulas of Theorem \ref{thm: Random formulas}. 

\begin{proof}[Proof of Theorem \ref{thm: Random formulas}]
	We begin with the first formula, which concerns the expectation of the number of tiles in random patches. Consider $G_{\bf \Sigma}$ equipped with weights as in Lemma \ref{lem: S,V,H,Q for random}, that is, if $\varepsilon=\varepsilon_T$ is an edge corresponding to a substitution tile $T\in\omega_k(T_i)$, then $w(\varepsilon_T)=p_{i,k}$. With weights assigned this way, the weight of a path in $G_{\bf \Sigma}$ is exactly the product of the probabilities for choosing the substitution rules associated with the edges it contains, and so applying Theorem \ref{thm: General graph counting} would give the required expectation formula: 
	
	First, assume that $S$ is a copy of a prototile $T_i\in\tau_{\bf \Sigma}$. For $0<a\le b\le 1$, a tile of type $r$ and volume in $[a,b]$ in a patch of the form $F^{\bf \Sigma}_t(S)$ corresponds to a walk of length $d\cdot t$ that initiates at vertex $i$ in $G_{\bf \Sigma}$ and terminates at a point on an interval contained in an edge $\varepsilon_T$, which is associated with a substitution tile $T$ of type $r$ in $\omega_k(T_h)$ for some $T_h\in\tau_{\bf \Sigma}$ and $1\le k\le \ell_h$. The edge $\varepsilon_T$ is of length $\log\frac{1}{\vol T}$ and terminates at vertex $r$. The distances of the boundaries of this interval from the initial vertex of $\varepsilon_T$ are given by
	\begin{align*}\label{eq: delta and eta}
	\delta_{\varepsilon_T}&=
	\begin{cases}
	\log\frac{1}{\vol T}-\log\frac{1}{a},&\vol T<a\le 1  \\
	0,&\rm{otherwise} 
	\end{cases}
	=\max\left(0, \log\tfrac{1}{\vol T}-\log \tfrac{1}{a}\right) \\
	\eta_{\varepsilon_T}&=
	\begin{cases}
	\log\frac{1}{\vol T}-\log\frac{1}{b},&\vol T<b\le  1\\
	0,&\rm{otherwise} 
	\end{cases}
	=\max\left(0, \log\tfrac{1}{\vol T}-\log \tfrac{1}{b}\right).\notag
	\end{align*}

	Theorem \ref{thm: General graph counting} now implies that
	\begin{align*}
	\frac{\mathbb{E}\left[\#\{\text{Type } r \text{ and }\vol\in[a,b]\text{ in } F_t^{\bf \Sigma}(S)\}\right]}{\vol(e^tS)}
	=\sum_{h=1}^{n}q^{\bf \Sigma}_{h}\sum_{k=1}^{\ell_h} p_{h,k}\sum\limits_{\substack{T\in\omega_k(T_h) \\ T\text { of type }r}}
	\int_{\delta_{\epsilon_T}}^{\eta_{\epsilon_T}}e^{-x}dx+o(1)
	\end{align*}
	as $t\rightarrow\infty$, with $q^{\bf \Sigma}_h$ any entry of column $h$ of $Q_{\bf \Sigma}$, and	where we use the fact that $\vol(e^tS)=e^{d\cdot t}$. Recalling the definition of $c_{\vol T}(x)$ in \eqref{eq: density functions} and substituting $u=\vol Te^x$ we get
	\begin{equation*}
	\int_{\delta_{\epsilon_T}}^{\eta_{\epsilon_T}}e^{-x}dx=\int_{\max(a,\vol T)}^{\max(b,\vol T)}\frac{\vol T}{u^2}du=\int_a^bc_{\vol T}(u)du.
	\end{equation*}
	Combined with Lemma \ref{lem: S,V,H,Q for random} and the definition of $C_{\bf \Sigma}(x)$ in \eqref{eq: density matrices} we conclude that
	\begin{align*}
	\frac{\mathbb{E}\left[\#\{\text{Type } r \text{ and }\vol\in[a,b]\text{ in } F_t^{\bf \Sigma}(S)\}\right]}{\vol(e^tS)}
	&=\int_a^b\frac{\left[{\bf v}^TC_{\bf \Sigma}(x)\right]_r}{{\bf v}^TH_{\bf \Sigma}\bf{1}}dx+o(1)
	\end{align*}
	as $t\rightarrow\infty$, for $S$ a copy of a prototile $T_i\in\tau_{\bf \Sigma}$, independent of $i$. Finally, if $S$ is a copy of $\alpha T_i$, then $F^{\bf \Sigma}_t(S)=F^{\bf \Sigma}_{t-\log(1/\vol S)}(T_i)$ and $e^{d\cdot t-\log(1/\vol S)}=\vol(e^tS)$. The first formula in Theorem \ref{thm: Random formulas} now follows for any $S$ of legal type and scale.
	
	For the second formula of Theorem \ref{thm: Random formulas}, which concerns the expectation of the volume occupied by tiles of a certain type, we extend ideas from the proof of the first formula and assign the following weights to the edges of $G_{\bf \Sigma}$: if $\varepsilon_T$ is an edge corresponding to a substitution tile $T\in\omega_k(T_i)$, then we set $w(\varepsilon_T)=p_{i,k}\cdot\vol T$.  With weights assigned in this way, the weight of a path in $G_{\bf \Sigma}$ is now the product of the probabilities for choosing the substitution rules associated with the edges it contains times the volume of the corresponding tile. 	With this choice of weights, the graph matrix function is 
	\begin{equation*}
	\left(M(s)\right)_{ij}=\sum\limits_{\substack{\varepsilon \text{ edge}\\i\rightarrow j}}w(\varepsilon)e^{-l(\varepsilon)s}=\sum_{k=1}^{\ell_i} p_{i,k}\sum\limits_{\substack{T\in\omega_k(T_i) \\ T\text { of type }j}}\vol T(\vol T)^{s}.
	\end{equation*}
	For every vertex $i$ in $G_{\bf\Sigma}$ the sum of the outgoing weights is exactly $1$, which implies that $\lambda=0$. Clearly $M(\lambda)=M(0)=M_{\bf \Sigma}(1)$ and $M'(0)=M'_{\bf \Sigma}(1)$ for $M_{\bf \Sigma}(s)$ as in Lemma \ref{lem: S,V,H,Q for random}, and so the matrix $Q$ associated with this choice of weights is equal to $Q_{\bf \Sigma}$. 
	
	For $S$ a copy of a prototile $T_i\in\tau_{\bf \Sigma}$, Theorem \ref{thm: General graph counting} now implies that
	\begin{align*}
	\frac{\mathbb{E}\left[\vol\left(\bigcup\text{Type } r \text{ and }\vol\in[a,b]\text{ in } F_t^{\bf \Sigma}(S)\right)\right]}{\vol(e^tS)}
	=\sum_{h=1}^{n}q^{\bf \Sigma}_{h}\sum_{k=1}^{\ell_h} p_{h,k}\sum\limits_{\substack{T\in\omega_k(T_h) \\ T\text { of type }r}}\vol T\int_{\delta_{\varepsilon_T}}^{\eta_{\varepsilon_T}}dx+o(1)
	\end{align*}
	as $t\rightarrow\infty$. Similarly to the proof of the first formula, with $d_{\vol T}(x)$ as in \eqref{eq: density functions}, we have
	\begin{align*}
	\int_{\eta_{\varepsilon_T}}^{\delta_{\varepsilon_T}}dx=\int_{\max(a,\vol T)}^{\max(b,\vol T)}\frac{\vol T}{u}du=\int_a^bd_{\vol T}(u)du.
	\end{align*}
	 We thus deduce that
	\begin{align*}
		\frac{\mathbb{E}\left[\vol\left(\bigcup\text{Tiles of type } r \text{ and }\vol\in[a,b]\text{ in } F_t^{\bf \Sigma}(S)\right)\right]}{\vol(e^tS)}
&=\int_a^b\frac{\left[{\bf v}^TD_{\bf \Sigma}(x)\right]_r}{{\bf v}^TH_{\bf \Sigma}\bf{1}}dx+o(1)
	\end{align*}
	as $t\rightarrow\infty$, for $S$ a copy of a prototile $T_i\in \tau_{\bf \Sigma}$, independently of $i$. The same argument as appeared in the first part of the proof implies that this is also true for any $S$ of legal type and scale, and so the proof of the second formula is complete.
	
	The final two formulas follow from the first two by setting $a=\min_{T\in\omega_{\bf \Sigma}}\vol T$ and $b=1$. Indeed, for every individual $T\in\omega_{\bf \Sigma}$ we have
	\begin{align*}
	\int_a^bc_{\vol T}(x)dx&=\int_{\vol T}^1
	\frac{\vol T}{x^2}dx=1-\vol T \\
	\int_a^bd_{\vol T}(x)dx&=\int_{\vol T}^1
	\frac{\vol T}{x}dx=-\vol T\cdot \log\vol T,
	\end{align*}
	and the formulas as stated follow from the definitions in \eqref{eq: S,V,H random}. 
\end{proof}

\appendix \section{The left Perron--Frobenius eigenvector ${\bf v}$ of $V_\sigma$} \label{sec: the role of v}

In the formulas presented in this paper, the entries of the left Perron--Frobenius eigenvector ${\bf v}\in\R^n$ are shown to determine the relative contribution of each of the $n$ different types of tiles to the various participating averages. We describe an additional geometric interpretation of ${\bf v}$ in terms of the {\it generation sequence of partitions} $\left(\delta_k\right)_{k\ge0}$ of $T_i\in\tau_\sigma$, generated by $\sigma$ according to a different procedure than that of Kakutani sequences. We restrict here to the non-random case to simplify notations.  

Set $\delta_0:=T_i$, and define $\delta_{k+1}$ by substituting {\bf all tiles} in $\delta_k$ simultaneously according to $\sigma$. The standard construction of substitution tilings, which includes the well-known Penrose and pinwheel tilings, is related to such sequences, see \cite{Baake-Grimm} for a comprehensive discussion. In the case of primitive (not necessarily normalized) fixed scale schemes, in which all the participating scales in \eqref{eq: substitution tiles} are identical, it is well-known that the Perron--Frobenius theorem implies the following formulas for the tile frequencies, as explained also in \cite[Theorem 6.12]{Yotam Kakutani}. These may be compared with the analogous formulas for incommensurable Kakutani sequences of partitions in Corollary \ref{cor: random kakutani frequencies}.

\begin{prop}\label{prop: fixed scale generation frequencies}
	Let $\xi$ be a primitive (not necessarily normalized) fixed scale substitution scheme in $\R^d$ with contraction constant $\alpha>0$, and let $\sigma$ be the equivalent normalized scheme. Let $\left(\delta_k\right)_{k\ge0}$ be a generation sequence of $T_i\in\tau_\sigma$ and let $1\le r \le n$. Then
	\begin{align}
	\#\{\emph{Tiles of type } r \emph{ in } \delta_k\}&=\frac{\vol T_i\cdot s_r}{{\bf s}^T{\bf u}_\vol}\alpha^{dk}+o\left(\alpha^{dk}\right),\notag\\
	\frac{\#\{\emph{Tiles of type } r \emph{ in } \delta_k\}}{\#\{\emph{Tiles in } \delta_k\}}&=s_r+o(1),\notag\\
	\vol\left(\bigcup\emph{Tiles of type } r \emph{ in } \delta_k\right)&=\vol T_i \cdot v_r+o(1),\notag
	\end{align}
	where ${\bf s},{\bf v}\in\R^n$ are left Perron--Frobenius eigenvectors of $S_\sigma$ and $V_\sigma$, respectively, normalized so that $\sum s_r=\sum v_r=1$, and ${\bf u}_\vol$ is the vector of volumes of prototiles. 
\end{prop}
In particular, since the third formula in Proposition \ref{prop: fixed scale generation frequencies} holds for any substitution scheme, the left Perron--Frobenius eigenvector ${\bf v}$ of $V_\sigma$, if chosen so that $\sum v_r=1$, can be viewed as the vector that registers the asymptotic volumes of the regions occupied by tiles of types $1\le r\le n$ in a {\bf generation} sequence of partitions generated by $\sigma$.

\subsection*{Acknowledgments}
I am grateful to Avner Kiro, Zemer Kosloff,  Joel Lebowitz, Uzy Smilansky, Yaar Solomon, Barak Weiss and Aron Wennman for helpful suggestions and discussions, and to the anonymous referee for valuable comments and presentational advice.

\end{document}